\documentclass[a4paper,12pt]{article}
\usepackage[english,russian]{babel}
\usepackage[cp1251]{inputenc}
\usepackage[dvipdfm]{graphicx}
\usepackage{mathtext}
\usepackage[T2A]{fontenc}
\usepackage{inputenc,amssymb}
\usepackage{epsfig}
\usepackage{caption2}
\usepackage{amsfonts}
\usepackage{amssymb}
\usepackage{amsmath}
\usepackage{geometry}
\usepackage{color}

\newcounter{eqcounter}[section]

\newcommand{\eqd}{\stackrel{d}{=}}

\newcommand{\R}{\mathbb R}
\newcommand{\I}{\mathbb{I}}

\title{Bounds for the accuracy of invalid normal approximation\footnote{Research supported by Russian Science Foundation, project 18-11-00155.}}

\author{Alexandra Dorofeeva\thanks{Faculty of Computational Mathematics and Cybernetics, Moscow State University, Moscow, Russia. E-mail: alex.dorofeyeva@gmail.com}, Victor Korolev\thanks{Faculty of Computational Mathematics and Cybernetics, Moscow State University, Moscow, Russia; Hanghzhou Dianzi University, Hangzhou, China; Federal Research Center ``Informatics and Conntrol'', Russian Academy of Sciences, Moscow, Russia. E-mail: vkorolev@cs.msu.ru}, Alexander Zeifman\thanks{Vologda State University, Vologda, Russia; Institute of Informatics Problems, Federal Research Center <<Computer Science and Control>> of the Russian Academy of Sciences, Moscow, Russia. E-mail: a$\_$zeifman@mail.ru}}

\date{}

\begin{document}

\maketitle

\noindent{\bf Abstract:} In applied probability, the normal approximation is often used for the distribution of data with assumed additive structure. This tradition is based on the central limit theorem for sums of (independent) random variables. However, it is practically impossible to check the conditions providing the validity of the central limit theorem when the observed sample size is limited. Therefore it is very important to know what the real accuracy of the normal approximation is in the cases where it is used despite it is theoretically inapplicable. Moreover, in some situations related with computer simulation, if the distributions of separate summands in the sum belong to the domain of attraction of a stable law with characteristic exponent less than two, then the observed distance between the distribution of the normalized sum and the normal law first decreases as the number of summands grows and begins to increase only when the number of summands becomes large enough. In the present paper an attempt is undertaken to give some theoretical explanation to this effect.

\smallskip

\noindent{\bf Key words:} central limit theorem; accuracy of the normal approximation; heavy tails; uniform distance; stable distribution

\smallskip

\noindent{\bf AMS 2000 subject classification:} 60F05, 60G50, 60G55, 62E20, 62G30

\section{Introduction}

In applied studies, the normal approximation is often used for the distribution of data with (at least assumed) additive structure. This tradition is based on the central limit theorem of probability theory which states that the distributions of sums of (independent) random variables satisfying certain conditions (say, the Lindeberg condition) converge to the normal law as the number of summands infinitely increases. However, it is practically impossible to check the conditions providing the validity of the central limit theorem when the observed sample size is limited. In particular, with moderate sample size, the histogram constructed from the sample from the Cauchy distribution whose tails are so heavy that even the mathematical expectation does not exist, is practically visually indistinguishable from the normal (Gaussian) density. Therefore it is very important to know what the real accuracy of the normal approximation is in the cases where it is used despite it is theoretically inapplicable. Moreover, in some situations related with computer simulation, if the distributions of separate summands in the sum belong to the domain of attraction of a stable law with characteristic exponent less than two, then the observed distance between the distribution of the normalized sum and the normal law first decreases as the number of summands grows and begins to increase only when the number of summands becomes large enough. In the present paper an attempt is undertaken to give some theoretical explanation to this effect. In Section 2 we introduce the notation, give necessary definitions and formulate some auxiliary results. In Section 3 the theorem is proved presenting the upper bound for the accuracy of the invalid normal approximation. In Section 4 the problem of evaluation of the threshold number of summands providing best possible accuracy of the invalid normal approximation is considered.

\section{Notation, definitions and auxiliary results}

Throughout the paper we assume that all the random variables are defined on the same probability space $(\Omega,\mathfrak{F},{\sf P})$. The mathematical expectation and variance with respect to the probability measure ${\sf P}$ will be denoted ${\sf E}$ and ${\sf D}$, respectively. The symbol $\eqd$ means the coincidence of distributions.

For $n\in\mathbb{N}$, let $X_1,\ldots,X_n$ be a homogeneous sample, that is, a set of independent identically distributed random variables with common distribution function $F(x)={\sf P}(X_1<x)$, $x\in\R$. For simplicity, without serious loss of generality we will assume that $F(x)$ is continuous.

We will follow the lines of approach described in \cite{Korolev2020}. Denote $S_n=X_1+\ldots+X_n$. The indicator of a set (event) $A\in\mathfrak{F}$ will be denoted $\I_A=\I_A(\omega)$, $\omega\in\Omega$:
$$
\I_A(\omega)=\begin{cases}1, & \omega\in A,\vspace{1mm}\cr 0, &\omega\notin A.\end{cases}
$$
Consider $u>0$ such that $0<F(u)<1$. It is obvious that $X_j=X_j\I_{\{|X_j|\le u\}}+X_j\I_{\{|X_j|> u\}}$. Then
$$
S_n=\sum\nolimits_{j=1}^n X_j\I_{\{|X_j|\le u\}}+\sum\nolimits_{j=1}^n X_j\I_{\{|X_j|> u\}}\equiv
S_n^{(\le u)}+S_n^{(> u)}.
$$
The number $N_n(u)$ of non-zero summands in the sum $S_n^{(\le u)}$ is a random variable that has the binomial distribution with parameters $n$ (``number of trials'') and $p=p(u)={\sf P}(|X_1|\le u)=F(u)-F(-u)$ (probability of ``success''). Note that, as $u$ infinitely grows, the parameter $p$ tends to 1. So, for $x\in\R$ we can write
$$
{\sf P}(S_n^{(\le u)}<x)\eqd {\sf P}\Big(\sum\nolimits_{j=0}^{N_n(u)}X_j^{(\le
u)}<x\Big),\eqno(1)
$$
where the random variables $X_1^{(\le u)},X_2^{(\le u)},\ldots$ are independent and have one and the same distribution function
$$
F^{(\le u)}(x)\equiv {\sf P}(X_1^{(\le u)}<x)={\sf
P}\big(X_1\I_{\{|X_1|\le u\}}<x\big|\,|X_1|\le u\big)=
$$
$$
=\frac{{\sf P}(X_1<x;\,|X_1|\le u)}{{\sf P}(|X_1|\le u)}=\begin{cases}1, &
x> u;\vspace{1mm}\\{\displaystyle\frac{F(x)-F(-u)}{F(u)-F(-u)}}, & |x|\le u; \vspace{1mm}\\0, & x<-u.\end{cases}\eqno(2)
$$
Moreover, the random variable $N_n(u)$ can be assumed to be independent of the sequence $X_1^{(\le u)},X_2^{(\le u)},\ldots$ For definiteness, if $N_n(u)=0$, then the sum $S_n^{(\le u)}$ is set equal to zero.

Similarly, for $x\in\R$ we have
$$
{\sf P}(S_n^{(>u)}<x)\eqd {\sf P}\Big(\sum\nolimits_{j=0}^{n-N_n(u)}X_j^{(>u)}<x\Big),\eqno(3)
$$
where $N_n(u)$ is {\it the same as in} (1) and is independent of the independent random variables $X_1^{(>u)},X_2^{(>u)},\ldots$ that have one and the same distribution function
$$
F^{(>u)}(x)\equiv {\sf P}(X_1^{(>u)}<x)={\sf P}\big(X_1\I_{\{|X_1|>u\}}<x\big|\,|X_1|>u\big)=
$$
$$
=\frac{{\sf P}(X_1<x;\,|X_1|> u)}{{\sf P}(|X_1|> u)}=\begin{cases}{\displaystyle\frac{F(x)}{F(-u)+1-F(u)}}, &
x<-u;\vspace{2mm}\\{\displaystyle\frac{F(-u)}{F(-u)+1-F(u)}}, & |x|\le u;\vspace{2mm}\\{\displaystyle\frac{F(-u)+F(x)-F(u)}{F(-u)+1-F(u)}}, & x>u.
\end{cases}\eqno(4)
$$
For definiteness, if $N_n(u)=n$, then the sum $S_n^{(>u)}$ is set equal to zero. Moreover, in (1) and (3) the random variables $X_1^{(\le u)},X_2^{(\le u)},\ldots,X_1^{(>u)},X_2^{(>u)},\ldots$ can be assumed to be jointly independent while the random variables $S_n^{(\le u)}$ and $S_n^{(> u)}$ are {\it not} independent and are related by the random variable $N_n(u)$.

\smallskip

{\sc Lemma 1.} {\it Let $A\in\mathfrak{F}$, $B\in\mathfrak{F}$. Then ${\sf P}(AB)\ge{\sf P}(A)-{\sf P}(\overline{B})$.}

\smallskip

The {\sc proof} is elementary.

\smallskip

The uniform (Kolmogorov) distance between the distribution functions $F_{\xi}$ and $F_{\eta}$ of random variables $\xi$ and $\eta$ will be denoted $\rho(F_{\xi},\,F_{\eta})$, $\rho(F_{\xi},\,F_{\eta})=\sup_x|F_{\xi}(x)-F_{\eta}(x)|$. The normal distribution function with expectation $a\in\R$ and variance $\sigma^2>0$ will be denoted $\Phi_{a,\sigma}(x)$,
$$
\Phi_{a,\sigma}(x)=\frac{1}{\sigma\sqrt{2\pi}}\int_{-\infty}^{x}\exp\Big\{-\frac{(z-a)^2}{2\sigma^2}\Big\}dz=\Phi_{0,1}\Big(\frac{x-a}{\sigma}\Big)=
\Phi_{0,\sigma}(x-a),\ \ \ x\in\R.
$$

\smallskip

{\sc Lemma 2.} {\it For any $a\in\R$, $\sigma>0$, $b\in\R$}
$$
\rho(\Phi_{a+b,\,\sigma},\,\Phi_{a,\,\sigma})= 2\Phi_{0,\sigma}\big({\textstyle\frac{|b|}{2}}\big)-1.
$$

\smallskip

{\sc Proof.} First, note that if $H(x)$ and $G(x)$ are two differentiable distribution functions, then $\rho(H,G)$ is realized (the supremum in $\sup_x|H(x)-G(x)|$ is attained) at one of the points $x$ where $F'(x)=G'(x)$. Indeed, we have
$$
\rho(H,G)=\sup_x|H(x)-G(x)|=\max\big\{\max_x\big[H(x)-G(x)\big], \max_x\big[G(x)-H(x)\big]\big\},
$$
and the extremum of each of the expressions in braces on the right-hand side is attained at the point where the derivative of the corresponding expression is equal to zero, which is equivalent to the equality of the derivatives of the distribution functions $H$ and $G$, that is, to the equality of the corresponding densities. In the case under consideration the latter condition is equivalent to that
$$
\frac{1}{\sigma\sqrt{2\pi}}\exp\Big\{-\frac12\Big(\frac{x-a-b}{\sigma}\Big)^2\Big\}=
\frac{1}{\sigma\sqrt{2\pi}}\exp\Big\{-\frac12\Big(\frac{x-a}{\sigma}\Big)^2\Big\},
$$
or $\big(x-(a+b)\big)^2=(x-a)^2$. Solving this equation we obtain $x-a=\frac{b}{2}$ yielding the desired result with the account of the relation $\Phi_{0,\sigma}(-|b|)=1-\Phi_{0,\sigma}(|b|)$.

\smallskip

Using the Lagrange formula, it is easy to deduce from Lemma 2 that
$$
\rho(\Phi_{a+b,\,\sigma},\,\Phi_{a,\,\sigma})\le \frac{|b|}{\sigma\sqrt{2\pi}}
$$
(see, e. g., inequality (3.4) in \cite{Petrov1972}).

\smallskip

{\sc Lemma 3}. {\it For $n\in\mathbb{N}$ let $\xi_1,\ldots,\xi_n$ be random variables, $a_1,\ldots,a_n$ be positive numbers such that
$a_1+\ldots+a_n=1$. Then for any $x>0$
$$
{\sf P}\Big(\Big|\sum\nolimits_{j=1}^n\xi_j\Big|\ge x\Big)\le\sum\nolimits_{j=1}^n{\sf P}(|\xi_j|\ge a_jx).
$$
If, in addition, the random variables $\xi_1,\ldots,\xi_n$ are identically distributed, then
$$
{\sf P}\Big(\Big|\sum\nolimits_{j=1}^n\xi_j\Big|\ge x\Big)\le n{\sf P}\big(|\xi_1|\ge {\textstyle\frac{x}{n}}\big).
$$
}

\smallskip

{\sc Proof}. First, note that
$$
{\sf P}\Big(\Big|\sum\nolimits_{j=1}^n\xi_j\Big|\ge x\Big)\le{\sf P}\Big(\sum\nolimits_{j=1}^n|\xi_j|\ge x\Big).
$$
Next, from geometrical considerations it follows that
$$
\Big\{\omega:\, \sum\nolimits_{j=1}^n|\xi_j(\omega)|\ge x\Big\}\subseteq\big\{\omega:\,|\xi_1(\omega)|\ge a_1x\big\}\bigcup
\Big\{\omega:\, \sum\nolimits_{j=2}^n|\xi_j(\omega)|\ge (1-a_1)x\Big\}\subseteq
$$
$$
\subseteq\big\{\omega:\,|\xi_1(\omega)|\ge a_1x\big\}\bigcup\big\{\omega:\,|\xi_2(\omega)|\ge a_2x\big\}\bigcup
\Big\{\omega:\, \sum\nolimits_{j=3}^n|\xi_j(\omega)|\ge (1-a_1-a_2)x\Big\}\subseteq\ldots
$$
$$
\ldots\subseteq\bigcup_{j=1}\nolimits^n\big\{\omega:\,|\xi_j(\omega)|\ge a_jx\big\}.
$$
Therefore,
$$
{\sf P}\Big(\Big|\sum\nolimits_{j=1}^n\xi_j\Big|\ge x\Big)\le{\sf P}\Big(\sum\nolimits_{j=1}^n|\xi_j|\ge x\Big)\le
$$
$$
\le{\sf P}\Big(\bigcup\nolimits_{j=1}^n\big\{\omega:\,|\xi_j(\omega)|\ge a_jx\big\}\Big)\le\sum\nolimits_{j=1}^n{\sf P}(|\xi_j|\ge a_jx).
$$
The lemma is proved.

\smallskip

{\sc Lemma 4}. {\it For $n\in\mathbb{N}$ let $\xi_1,\ldots,\xi_n$ be random variables such that ${\sf E}|\xi_j|^{\delta}<\infty$ for some $\delta>0$, $j=1,\ldots,n$. Denote $\theta_n=\xi_1+\ldots+\xi_n$.

\noindent {\rm(i)} If $0<\delta\le1$, then
$$
{\sf E}|\theta_n|^{\delta}\le \sum\nolimits_{j=1}^n{\sf E}|\xi_j|^{\delta}.
$$

\noindent {\rm (ii)} If $1\le\delta\le2$, the random variables $\xi_1,\ldots,\xi_n$ are independent and ${\sf E}\xi_j=0$, $j=1,\ldots,n$, then
$$
{\sf E}|\theta_n|^{\delta}\le \Big(2-\frac1n\Big)\sum\nolimits_{j=1}^n{\sf E}|\xi_j|^{\delta}.
$$
}

\smallskip

{\sc Proof}. Statement (i) is elementary, statement (ii) was proved in \cite{BahrEsseen1965}.

\section{Main results}

Consider the upper bound for the uniform distance between the distribution of the normalized sum
$$
S_n^*=\frac{1}{\sqrt{n}}\sum\nolimits_{j=1}^nX_j
$$
and the normal law with some expectation $a\in\R$ and variance $\sigma^2>0$. The choice of concrete values of $a$ and $\sigma^2$ will be discussed later.

From what has been said it follows that
$$
S_n^*\eqd\frac{S_n^{(\le u)}}{\sqrt{n}}+\frac{S_n^{(> u)}}{\sqrt{n}}.
$$
For brevity and convenience, we will use the notation
$$
\zeta_n=\frac{S_n^{(\le u)}}{\sqrt{n}},\ \ \ \eta_n=\frac{S_n^{(> u)}}{\sqrt{n}}.
$$

\smallskip

{\sc Theorem 1.} {\it Let $u>0$ be arbitrary. Then for any $a\in\R$ and $\sigma>0$ we have}
$$
\rho(F_{\zeta_n+\eta_n},\,\Phi_{a,\,\sigma})\le\rho(F_{\zeta_n},\,\Phi_{a,\,\sigma})+n\big(F(-u)+1-F(u)\big).\eqno(5)
$$

\smallskip

{\sc Proof}. Let $\epsilon>0$ be arbitrary. According to Lemma 1 we have
$$
{\sf P}(\zeta_n+\eta_n<x)={\sf P}(\zeta_n+\eta_n<x;\,|\eta_n|\le\epsilon)+{\sf P}(\zeta_n+\eta_n<x;\,|\eta_n|>\epsilon)\ge
$$
$$
\ge{\sf P}(\zeta_n<x-\eta_n;\,|\eta_n|\le\epsilon)\ge{\sf P}(\zeta_n<x-\epsilon;\,|\eta_n|\le\epsilon)\ge{\sf P}(\zeta_n<x-\epsilon)-{\sf P}(|\eta_n|\ge\epsilon).\eqno(6)
$$
On the other hand, obviously,
$$
{\sf P}(\zeta_n+\eta_n<x)={\sf P}(\zeta_n<x-\eta_n;\,|\eta_n|\le\epsilon)+{\sf P}(\zeta_n+\eta_n<x;\,|\eta_n|>\epsilon)\le
$$
$$
\le{\sf P}(\zeta_n<x+\epsilon;\,|\eta_n|\le\epsilon)+{\sf P}(\zeta_n+\eta_n<x;\,|\eta_n|>\epsilon)\le{\sf P}(\zeta_n<x+\epsilon)+{\sf P}(|\eta_n|>\epsilon).\eqno(7)
$$
It is easy to see that
$$
|{\sf P}(\zeta_n+\eta_n<x)-\Phi_{a,\,\sigma}(x)|=\max\big\{{\sf P}(\zeta_n+\eta_n<x)-\Phi_{a,\,\sigma}(x),\,\Phi_{a,\,\sigma}(x)-{\sf P}(\zeta_n+\eta_n<x)\big\}.\eqno(8)
$$
Using (7) and Lemma 2 we obtain
$$
{\sf P}(\zeta_n+\eta_n<x)-\Phi_{a,\,\sigma}(x)\le {\sf P}(|\eta_n|>\epsilon)+\big[{\sf P}(\zeta_n<x+\epsilon)-\Phi_{a,\,\sigma}(x+\epsilon)\big]+
$$
$$
+\big[\Phi_{a,\,\sigma}(x+\epsilon)-\Phi_{a,\,\sigma}(x)\big]\le{\sf P}(|\eta_n|>\epsilon)+\rho(F_{\zeta_n},\,\Phi_{a,\,\sigma})+
\big[2\Phi_{0,\sigma}\big({\textstyle\frac{\epsilon}{2}}\big)-1\big].\eqno(9)
$$
Using (6) and Lemma 2 we obtain
$$
\Phi_{a,\,\sigma}(x)-{\sf P}(\zeta_n+\eta_n<x)\le \Phi_{a,\,\sigma}(x)-{\sf P}(\zeta_n<x-\epsilon)+{\sf P}(|\eta_n|>\epsilon)=
$$
$$
=\big[\Phi_{a,\,\sigma}(x)-\Phi_{a,\,\sigma}(x-\epsilon)\big]-\big[{\sf P}(\zeta_n<x-\epsilon)-\Phi_{a,\,\sigma}(x-\epsilon)\big]+ {\sf P}(|\eta_n|>\epsilon) \le
$$
$$
\le\rho(F_{\zeta_n},\,\Phi_{a,\,\sigma})+{\sf P}(|\eta_n|>\epsilon)+\big[2\Phi_{0,\sigma}\big({\textstyle\frac{\epsilon}{2}}\big)-1\big].\eqno(10)
$$
Substituting (9) and (10) in (8) we obtain
$$
\rho(F_{\zeta_n+\eta_n},\,\Phi_{a,\,\sigma})\le\rho(F_{\zeta_n},\,\Phi_{a,\,\sigma})+\big[2\Phi_{0,\sigma}\big({\textstyle\frac{\epsilon}{2}}\big)-1\big]+
{\sf P}(|\eta_n|>\epsilon).\eqno(11)
$$
To estimate the last term on the right-hand side of (11) use representation (3) for $S_n^{(>u)}$. With the account of the convention $\sum_{j=0}^0=0$, by the formula of total probability we have
$$
{\sf P}(|\eta_n|>\epsilon)={\sf P}\Big(\Big|\sum\nolimits_{j=1}^{n-N_n(u)}X_j^{(>u)}\Big|>\epsilon\sqrt{n}\Big)=
$$
$$
=\sum_{k=1}^nC_n^k\big(F(-u)+1-F(u)\big)^k
\big(F(u)-F(-u)\big)^{n-k}{\sf P}\Big(\Big|\sum\nolimits_{j=1}^kX_j^{(>u)}\Big|>\epsilon\sqrt{n}\Big).\eqno(13)
$$
Estimating the probability in the last expression by Lemma 3, continuing (13) we obtain
$$
{\sf P}(|\eta_n|>\epsilon)\le\sum_{k=1}^nC_n^k\big(F(-u)+1-F(u)\big)^k
\big(F(u)-F(-u)\big)^{n-k}k{\sf P}\big(|X_j^{(>u)}|>{\textstyle\frac{\epsilon\sqrt{n}}{k}}\big)\le
$$
$$
\le{\sf P}\big(|X_j^{(>u)}|>{\textstyle\frac{\epsilon}{\sqrt{n}}}\big)\cdot\sum_{k=1}^nC_n^k\big(F(-u)+1-F(u)\big)^k
\big(F(u)-F(-u)\big)^{n-k}k=
$$
$$
=n\big(F(-u)+1-F(u)\big){\sf P}\big(|X_j^{(>u)}|>{\textstyle\frac{\epsilon}{\sqrt{n}}}\big).
$$
Substitution of this bound in (11) yields
$$
\rho(F_{\zeta_n+\eta_n},\,\Phi_{a,\,\sigma})\le\rho(F_{\zeta_n},\,\Phi_{a,\,\sigma})+\big[2\Phi_{0,\sigma}\big({\textstyle\frac{\epsilon}{2}}\big)-1\big]+
n\big(F(-u)+1-F(u)\big){\sf P}\big(|X_j^{(>u)}|>{\textstyle\frac{\epsilon}{\sqrt{n}}}\big).\eqno(14)
$$
Now let $\epsilon\to 0$ in (14) and obtain the desired result. The theorem is proved.

\smallskip

In practice, the values of the parameters $a$ and $\sigma$ can be chosen by the following reasoning. It is easy to verify (say, by the consideration of characteristic functions) that
$$
S_n^{(\le u)}\eqd\sum\nolimits_{j=1}^n\widetilde{X}_j^{(\le u)},
$$
where $\widetilde{X}_1^{(\le u)},\ldots,\widetilde{X}_n^{(\le u)}$ are independent identically distributed random variables,
$$
\widetilde{X}_j^{(\le u)}=\begin{cases}X_j^{(\le u)} & \text{ with probability } F(u)-F(-u);\vspace{2mm}\\ 0 & \text{ with probability } F(-u)+1-F(u).\end{cases}
$$
Then in accordance with (3), the parameter $a$ can be defined as
$$
a=a(u)={\sf E}\widetilde{X}_1^{(\le u)}=[F(u)-F(-u)]{\sf E}X_1^{(\le u)},
$$
and the parameter $\sigma^2$ can be defined as
$$
\sigma^2=\sigma^2(u)={\sf D}\widetilde{X}_1^{(\le u)}=[F(u)-F(-u)]{\sf D}X_1^{(\le u)}+[F(-u)+1-F(u)]\big({\sf E}X_1^{(\le u)}\big)^2.
$$
With these values of $a$ and $\sigma$ the first term on the right-hand side of (5) will tend to zero by the central limit theorem as $n\to\infty$, and can be estimated by the standard techniques, say, by the Berry--Esseen inequality for binomial random sums, see \cite{KorolevDorofeeva2017, KorolevShevtsova2012}.

As regards the second term on the right-hand side of (5), with large $u$, $p=F(u)-F(-u)$ close to one and moderate (but large enough) $n$ the term $\eta_n$ may be small due to that the sum $S_n^{(>u)}$ contains very few summands. Moreover, in the case of light tails, putting $u=u_n$ so that $n[F(-u)+1-F(u_n)]\to 0$ as $n\to\infty$, it is possible to make sure that the right-hand side of (5) can be made arbitrarily small by the choice of arbitrarily large $n$ so that the limit distribution for the normalized sum $S_n^*$ will be normal.

Under some additional conditions, at the expense of introducing additional parameter, the dependence of the second term of the bound given in Theorem 1 on $n$ can be made better.

\smallskip

For $c\in(0,2]$ let $h(c)=\I_{(1,2]}(c)$.

\smallskip

{\sc Theorem 2}. {\it Assume that the distribution function $F(x)$ belongs to the domain of attraction of a stable law with characteristic exponent $\alpha\in(0,2)$. If, moreover, $\alpha\ge1$, then additionally assume that $F$ is symmetric $($that is, $F(-x)=1-F(x)$ for $x>0$$)$. Then for any $u>0$ $\epsilon>0$ and $\delta\in(0,\alpha)$ we have}
$$
\rho(F_{\zeta_n+\eta_n},\,\Phi_{a,\,\sigma})\le\rho(F_{\zeta_n},\,\Phi_{a,\,\sigma})+\big[2\Phi_{0,\sigma}\big({\textstyle\frac{\epsilon}{2}}\big)-1\big]+
2^{h(\delta)}\epsilon^{-\delta}n^{1-\delta/2}\big(F(-u)+1-F(u)\big){\sf E}\big|X_1^{(>u)}\big|^{\delta}.\eqno(15)
$$

\smallskip

{\sc Proof}. The starting point of the proof is inequality (11). In accordance with the result of \cite{Tucker1975}, in the case under consideration ${\sf E}|X_1|^{\delta}<\infty$ for any $\delta\in(0,\alpha)$, and hence, ${\sf E}|X_1^{(>u)}|^{\delta}<\infty$. Moreover, if $\alpha>1$, then the mathematical expectation of $X_1$ exists and, due to the assumption that in that case the distribution of $X_1$ is symmetric, ${\sf E}X_1=0$. Therefore by representation (3), the Markov inequality, and Lemma 4, continuing (13) we obtain
$$
{\sf P}(|\eta_n|>\epsilon)\le\frac{1}{\epsilon^{\delta}n^{\delta/2}}\sum\nolimits_{k=0}^{n}C_n^k\big(F(-u)+1-F(u)\big)^k\big(F(u)-F(-u)\big)^{n-k}
{\sf E}\Big|\sum\nolimits_{j=0}^kX_j^{(>u)}\Big|^{\delta}\le
$$
$$
\le\frac{2^{h(\delta)}{\sf E}\big|X_1^{(>u)}\big|^{\delta}}{\epsilon^{\delta}n^{\delta/2}}\sum\nolimits_{k=0}^nC_n^kk\big(F(-u)+1-F(u)\big)^k\big(F(u)-F(-u)\big)^{n-k}=
$$
$$
=2^{h(\delta)}\epsilon^{-\delta}n^{1-\delta/2}\big(F(-u)+1-F(u)\big){\sf E}\big|X_1^{(>u)}\big|^{\delta}.\eqno(16)
$$
The theorem is proved.

\smallskip

We see that in (15) the exponent of $n$ is less than that in (5). However, in (15) an additional parameter $\epsilon$ appeared. The second term on the right-hand side of (15) can be made arbitrarily small by the appropriate choice of $\epsilon$. With $n$ and $\epsilon$ fixed, the third term on the right-hand side of (15) can be made arbitrarily small by the choice of large enough $u$.

Actually Theorems 1 and 2 are simple variants of a so-called pre-limit theorem, see \cite{Klebanov1999}.

As an illustration of how Theorem 2 acts, consider the following example.

\smallskip

{\sc Example}. Assume that the random variables $X_1,X_2,\ldots$ have common probability density
$$
f(x)=\frac{3}{4(1+|x|)^{5/2}},\ \ \ x\in\R.
$$
The corresponding distribution function has the form
$$
F(x)=\begin{cases}{\displaystyle\frac{1}{2(|x|+1)^{3/2}}}, & x<0;\vspace{2mm}\\
{\displaystyle 1-\frac{1}{2(x+1)^{3/2}}}, & x\ge0.\end{cases}
$$
This distribution function belongs to the domain of attraction of a stable law with characteristic exponent $\frac32$. It is easy to make sure that in this case, according to (4), we have
$$
{\sf E}|X_1^{(>u)}|=2\int_{0}^{\infty}xdF^{(>u)}(x)={\textstyle\frac32}(u+1)^{3/2}\int_u^{\infty}\frac{xdx}{(1+x)^{5/2}}=3(u+1)-1.
$$
Hence, choosing $\delta=1$ we see that the right-hand side of (16) is
$$
\frac{\sqrt{n}(1-p){\sf E}|X_1^{(>u)}|}{\epsilon}=\frac{\sqrt{n}(3u-2)}{2\epsilon(u+1)^{3/2}}=\frac{\sqrt{n}}{2\epsilon}\Big[\frac{3}{\sqrt{u+1}}+
o\Big(\frac{1}{\sqrt{u}}\Big)\Big]
$$
as $u\to\infty$. Therefore, with fixed $\epsilon>0$ and $n$ by choosing $u$ large enough the third summand on the right-hand side of (15) can be made arbitrarily small.

\section{On the threshold value of the number of summands}

Consider the problem of determination of $n_0$ such that for $n$ growing from 1 to $n_0$ the distance $\rho(F_{\zeta_n+\eta_n},\,\Phi_{a,\sigma})$ decreases and for $n>n_0$ this distance increases. Assume that the first summand on the right-hand side of (5) is estimated by the Berry--Esseen inequality
$$
\rho(F_{\zeta_n},\,\Phi_{a,\sigma})\le \frac{C_0\widetilde{L}_3^{(\le u)}}{\sqrt{n}},
$$
where $\widetilde{L}_3^{(\le u)}$ is the Lyapunov fraction,
$$
\widetilde{L}_3^{(\le u)}=\frac{{\sf E}\big|\widetilde{X}_1^{(\le u)}-{\sf E}\widetilde{X}_1^{(\le u)}\big|^3}{\big({\sf D}\widetilde{X}_1^{(\le u)}\big)^{3/2}},
$$
$C_0>0$ is the absolute constant, $C_0\le 0.4690$ \cite{Shevtsova2014}. It is easy to verify that if $c>0$, $d>0$, then
$$
\mathrm{arg}\min_{z>0}\Big(\frac{c}{\sqrt{z}}+dz\Big)=\Big(\frac{c}{2d}\Big)^{2/3}.
$$
Putting $z=n$, $c=C_0\widetilde{L}_3^{(\le u)}$, $d=F(-u)+1-F(u)$, we see that the minimum of the upper bound for $\rho(F_{\zeta_n+\eta_n},\,\Phi_{a,\sigma})$ is attained at $n_0$ which is either the integer part of
$$
z_0=\bigg(\frac{C_0\widetilde{L}_3^{(\le u)}}{2\big(F(-u)+1-F(u)\big)}\bigg)^{2/3},
$$
or at $n_0+1$.

Now assume that conditions of Theorem 2 hold. Then we have
$$
\mathrm{arg}\min_{z>0}\Big(\frac{c}{\sqrt{z}}+dz^{1-\delta/2}\Big)=\bigg(\frac{c}{2d(1-\frac{\delta}{2})}\bigg)^{\frac{2}{3-\delta}}
$$
so that
$$
z_0=\bigg(\frac{C_0\widetilde{L}_3^{(\le u)}}{2(1-\frac{\delta}{2})\big(F(-u)+1-F(u)\big)}\bigg)^{\frac{2}{3-\delta}},
$$

\renewcommand{\refname}{References}

\end{document}